\documentclass[twoside]{article}

\input{spsp12.sty}

\def\E {{\rm E}}

\def\GSN {{\rm GSN}}

\def\N {{\rm N}}
\def\HN {{\rm HN}}
\def\Var {{\rm Var}}
\def\Cov {{\rm Cov}}

\def\diag {{\rm diag}}

\newcommand{\bm}[1]{\mbox{\boldmath{$#1$}}}
\begin{document}
\title{\bf \vspace{2cm}
Exploring the Distributional Properties of the Non-Gaussian Random Field Models}
\date{}
\author{ Behzad Mahmoudian\footnote{Speaker: b.mahmoudian@qom.ac.ir} $^1$\\
\small $^1$Department of Statistics, Faculty of Science, University of Qom}
\fancyhead[LE] {Exploring the Distributional Properties of Shape-mixture of SGRFs}
\fancyhead[RO] {Mahmoudian, B.} 
\maketitle
\firstpage

\abstract In the environmental modeling field, the exploratory analysis of
responses often exhibits spatial correlation as well as some
non-Gaussian attributes such as skewness and/or heavy-tailedness.
Consequently, we
propose a general spatial model based on scale-shape mixtures of the
multivariate skew-normal distribution. Intuitively, it incorporates
distinct random effects to account for the spatial dependencies not
explained by a simple Gaussian random field model. Importantly, the
proposed model is capable of generating a wide range of skewness and kurtosis levels.
Meanwhile, we demonstrate that the skewness mixing can induce
asymmetric tail dependence at sub-asymptotic and/or asymptotic
levels.\\ \\
\keyword Scale-shape mixtures, Heavy tails, Skewness, Tail dependence.

\ams  62M30, 62G32.

\hspace{-1cm}\rule{\textwidth}{0.2mm}

\section{Introduction}
Many environmental random phenomena are spatially correlated, which
means that similar observed values in a domain are likely to occur
near one another than those far away. Random field (RF) models have
been widely investigated for analysis of this kind of data that
arise in epidemiology, climatology and many other disciplines. The
Gaussian random field (GRF) model is fairly well-accepted as a
custom working model, see
\citet{Gel1} for more detailed information. Despite such well-established theory, GRF models are
not always preferred in studies of empirical data that do not
conform to bell-shape distributions. In
other words, spatial responses usually exhibit substantial skewness
and/or extra kurtosis, which is particularly prevalent in
environmental applications.

One of standard and pragmatic approaches is to consider the
multivariate non-normal family of distributions to extend the GRF
model. In this regard, the specification of a distribution family
with particular behaviour for the finite-dimensional distributions
of a RF is a delicate issue, see \cite{Ma1} and \cite{Min}. Much progress has been made in the
general area of non-Gaussian RF models, such as elliptically
contoured RF models \citep{Ma2}, skew-Gaussian RF (SGRF)
models \citep{ZEL,Gent,Sch}, asymmetric Laplace RF models \citep{Lum}
and Tukey g-and-h RF models \citep{Xu2},
just to name a few recent contributions. These models are applied to
accommodate skewness and/or heavy-tailedness encountered in spatial
data. Even though these non-Gaussian models lead to desirable
modeling strategies, it is not guaranteed that they should always be
applicable. For example, their model structures can only induce
limited range of either but not both of positive or negative skewness.
Further, their model performances are hindered by constraints on the
parameters to ensure the existence of the moments.

When a spatial tactic is applied to extreme data, its attainable
degree of tail dependence must be characterized and quantified. This
is a big concern in dependence modeling
\citep{App,Dav2,Opt,Hus1,Mor,Wad2,Kru2,Hus2}. Max-stable RF models
are a useful tool to analyze spatial extremes and widely considered
to model the maxima observed at sites in a spatial domain, see
\cite{Dey} and the references therein for further details. Unfortunately,
these models assume that the marginal variables are asymptotically
dependent with dependence structure determined by rigid form taken
by the asymptotic results. Such an assumption is inappropriate for
real-world applications. Actually, fitting a misspecified model to
the data contributes to an incorrect estimation of probabilities of
extreme joint events. Practically speaking, the asymptotic arguments
as well as statistical inference for tail dependence analysis is
applicable as the number of independent replications from the
underlying field becomes large. On the other hand, some of these
models may be appropriate only on the local scale where observations
collected over a small number of spatial locations are assumed {\it
a priori} to be always dependent. Moreover, a general class of
models that induce asymmetric tail dependence at sub-asymptotic
and/or asymptotic levels have so far received little attention.
These issues motivate us to develop a variant of SGRF which
generates more sophisticated tail dependence structures for extremal
data or other heavy-tailed phenomena.

The current work is built on earlier study by \citet{Mah1}, who
employed a four-level hierarchical spatial model in terms of the
generalized skew-normal distribution (GSN) of \citet{Sah}, but here
the mixing components which incorporate the skewness is embedded in
the first stage of hierarchy, supporting plausible estimation
results. Besides, the ideas inspired by this author are here
extended to estimate the direction of skewness from data.
Because the accessible skewness under GSN distribution
is limited, we consider the scale-mixtures of this probability model to
induce an unlimited amount of the skewness. Outliers as well as
regions with inflated variances may be detected in the Gaussian
framework, by virtue of taking into account the scale mixtures of
GRF models \citep{Pal,Bu,Fag}. To model simultaneously
skewness and heaviness in tails, the GSN distributed RF model is
rescaled according to the suggested model of \citet{Pal} and is
reshaped in terms of a GRF model. The final step, i.e.
shape/skewness mixture formulation is employed to address the
challenge of skewness direction and magnitude identification in
spatial modeling.

We illustrate that our model is capable of describing
various amount of skewness and kurtosis ranging from mild to large.
We realize that not only all finite-dimensional distributions of the
proposed SGRF model are asymptotically independent but also at
finite levels different degrees of dependence are achievable.
Fortunately, the skewness parameters play the main role in this
respect so that the amount of skewness towards different directions
calibrates speed of convergence of the tail probability to the
asymptotically independent and/or dependent limits. We hope that
using the model with the aforementioned tail characteristics could
support sensible risk estimation of severe joint extreme events.
Moreover, the parametrization of the adopted distribution for a SGRF
is such that the second-order stationarity assumption is not
violated and its covariances vanish as the distances among spatial
locations go to infinity.

The article is organized as follows. In Section \ref{model1}, we introduce the GSN
distribution as a skew-normal model of interest.
Then in Section \ref{model2}, we discuss about the tail dependence properties, moments and stochastic
representation of the GSN distribution. The robustness of the proposed RF model is also
studied in Section \ref{model3}. Finally, conclusions based on the
results are given in Section \ref{conclu}.
\section{The GSN distribution}\label{model1}
One of the challenges for statistical procedures is to define skewed
distributions. In the large class of skew models
\citep[e.g.,][]{Az}, we restrict attention to the GSN family of
distributions. Reasons behind this are preservation of the
correlation structure under induced skewness, appealing generating
mechanisms and desirable fitting properties. Let
$\phi_n(\cdot;\bm\mu,\bm\Sigma)$
and $\Phi_n(\cdot;\bm\mu,\bm\Sigma)$ are the probability distribution function
(pdf) and cumulative
distribution function (cdf) of $\N_n(\bm\mu,\bm\Sigma)$,
respectively. Concerning its
definition, a n-dimensional random vector $\bm Z$ is said to have a
multivariate GSN distribution, denoted by $\GSN_n(\bm
\mu,\bm\Sigma,\bm\delta)$, if its pdf is of the form
\begin{align}\label{SSN}
f(\bm z)=2^n\phi_n(\bm z;\bm\mu,\bm\Sigma+\bm D^2)\Phi_n(\bm
D(\bm\Sigma+\bm D^2)^{-1}(\bm z-\bm\mu);\bm0,\bm\Delta),\quad \bm
z\in \mathbb{R}^n,
\end{align}
where $\bm\mu\in \mathbb{R}^n$, $\bm\Sigma\in \mathbb{R}^{n\times
n}$, $\bm\delta\in \mathbb{R}^n$, $\bm D=\diag(\bm\delta)$, $\bm
D^2=\diag(\delta_1^2,\ldots,\delta_n^2)$ and $\bm\Delta=\bm I_n-\bm
D(\bm\Sigma+\bm D^2)^{-1}\bm D$. Here, $\diag(\bm\eta)$ represents a diagonal matrix
with diagonal elements specified by the vector $\bm\eta$.
Note that for $\bm\delta=\bm 0_n$ where $\bm 0_n$ is a
$n\times 1$ vector of zeros, \eqref{SSN} reduces to the symmetric
$\N_n(\bm\mu,\bm\Sigma)$ pdf, whereas for non-zero values of
$\bm\delta$, it produces a perturbed family of
$\N_n(\bm\mu,\bm\Sigma)$ pdfs. If $\bm Z$ has pdf \eqref{SSN},
its moment generating function (mgf) is given implicitly by
\begin{align}
M_Z(\bm t)=2^n\exp\left(\bm t'\bm\mu+\frac{1}{2}\bm t'(\bm\Sigma+\bm
D^2)\bm t\right)\Phi_n(\bm D\bm t),\quad \bm t\in
\mathbb{R}^n,\label{mgfssn}
\end{align}
in which $\Phi_n(\cdot)$ is the cdf of $\N_n(\bm 0,\bm I_n)$. Let
$V_i\stackrel{\;\rm i.i.d.}{\sim}\HN_1(0,1)$
where $\HN_1$ represents univariate half standard normal distribution
and $\bm V=(V_1,\ldots,V_n)'$ be independent of $\bm W\sim\N_n(\bm 0,\bm\Sigma)$.
The GSN distribution as defined in \eqref{SSN} would
be stochastically represented as $\bm Z\stackrel{\,\rm d}{=}\bm\mu+\bm D\bm V+\bm W$
in which $\stackrel{\,\rm d}{=}$ means `as distributed'.

In the remaining part of this section, we discuss about the tail
probabilities of the GSN model. The tail dependence coefficient
(TDC) is a simple measure to quantify occurrences of the concurrent
extreme events; high level of TDC implies more probability of
simultaneous extreme events. We focus on the
upper tail of the GSN distribution; the lower tail properties can be
considered similarly. Let $\bm Z=(Z_1,Z_2)'$ be a two-dimensional
random vector. The upper tail dependence coefficient (UTDC) of a
random vector $\bm Z$ is defined by
\begin{align*}
\chi=\lim_{u\rightarrow 1^-}P[F_2(Z_2)>u|F_1(Z_1)>u],
\end{align*}
where $F_i(\cdot)$ for $i=1,2$ is the marginal cdf of $Z_i$.
The bivariate distribution family is said to be upper
tail dependent if $0<\chi\leq 1$ and upper tail independent if
$\chi=0$, in the case the limit exists. In particular, the
multivariate normal distribution cannot accommodate tail dependency
\citep{Col}. From \citet{Ber2} as well as the
references therein, we know that most of the multivariate
skew-normal distributions are asymptotically independent. Unfortunately,
the default version of the GSN distribution entirely lacks any flexibility in tail
dependence. Thereby, we investigate the tail probabilities of the
following bivariate GSN distribution
\begin{align}\label{biGSN}
\GSN_2\left(-\sqrt{\frac{2}{\pi}}\bm\Gamma^{1/2}\bm\delta, \bm\Gamma,\bm\Gamma^{1/2}\bm\delta
\right),
\end{align}
where $\bm\Gamma$ is an $2\times2$ correlation matrix, whose
off-diagonal elements are equal to $\rho$. 
Since the tail dependence only depends on the tail behaviour of the random
variables, the GSN distributed random vector with zero mean vector
is designated in \eqref{biGSN}. The intuition behind
this specification is that the tail flexibility becomes possible in
some particular setting when
a shape-mixture extension of the GSN distribution is taken into consideration.
We study the tail property for the skew-normal model of interest in the following proposition.

{\pro\label{pro1} The UTDC of the GSN distribution in \eqref{biGSN} is zero for
\begin{itemize}
\item[(a)]{$0\leq\delta_1\leq\delta_2$, $\delta_1,\delta_2<0$ and $\delta_1<0\leq\delta_2$.}
\item[(b)]{$0\leq\delta_2<\delta_1$ whenever
\begin{align*}
\delta_1<\sqrt{\frac{(1+\delta^2_2)(1+\rho)}{2\rho}-1}.
\end{align*}
}
\end{itemize}}
\textbf{Proof.} See \cite{Mah3} for a proof.\\
According to Proposition \ref{pro1} and its proof, the regular arguments do not entail the asymptotic independence for
$0\leq\delta_2<\delta_1$ with condition $\delta_1>\sqrt{(1+\delta^2_2)(1+\rho)/(2\rho)-1}$, and $\delta_2<0\leq\delta_1$.
\begin{figure*}[!t]
\vspace{0mm}
\begin{center}
\includegraphics[scale=0.9]{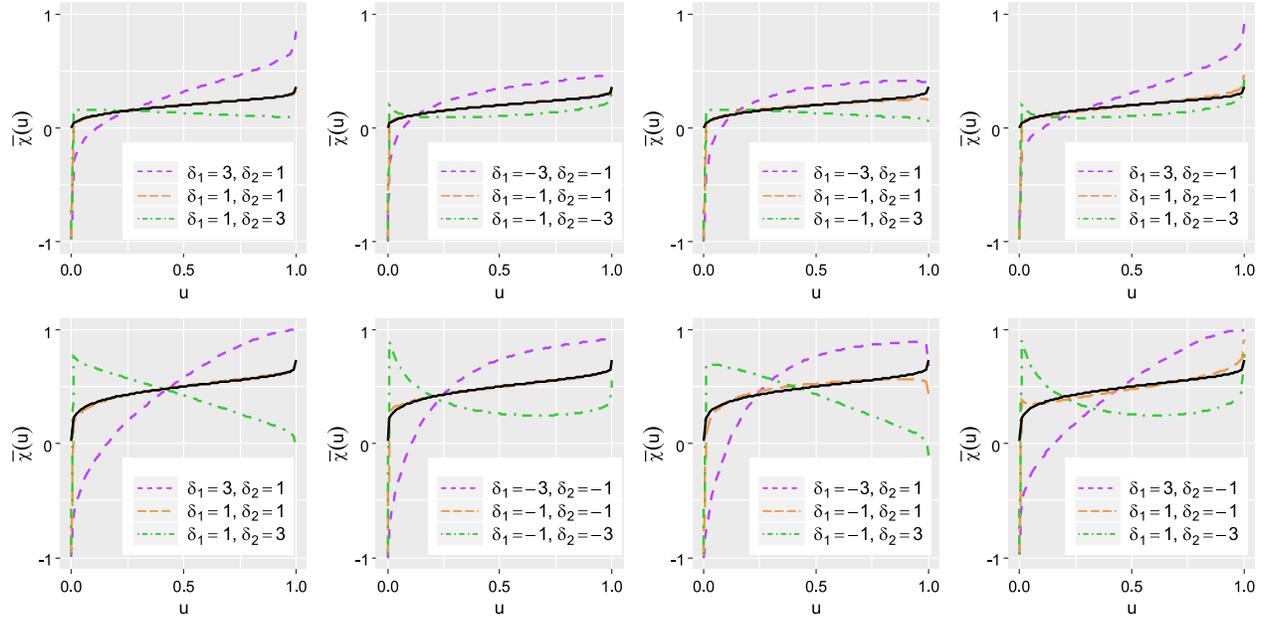}
 \end{center}
 \caption{The dependence measure $\bar{\chi}(u)$ for the GSN distribution: the curves shown on first row correspond to $\rho=0.4$
 as well as specified values of the $\delta_i$. The second row accords with $\rho=0.8$.  The solid line in each panel corresponds to
 the normal distribution.}\label{utdGSN1}
\end{figure*}

Under the asymptotic independence, \citet{Col} recommended to
characterize the extremal dependence at finite upper levels by
\begin{align*}
\bar{\chi}(u)=\frac{2\ln(P[F_1(Z_1)>u])}{\ln(P[F_1(Z_1)>u,F_2(Z_2)>u])}-1,
\end{align*}
where $-1\leq \bar{\chi}(u)\leq 1$ for all $0\leq u\leq1$ and its
larger magnitudes correspond to stronger dependence.
\eqref{utdGSN1} illustrate how the tail
probabilities under the given GSN distribution in \eqref{biGSN}
varies with different values of the correlation parameter $\rho$
and various choices of the skewness
parameters $\delta_i$.

These figures reveal a large variety of tail behaviours with
different decay or increase rates. One can see that the sign change
of skewness parameters results in conflicting shapes of the curves such as
position of the curves relative to the Gaussian counterpart and
dense or sparse behaviour. Another remark is that the tail
probabilities of the GSN distribution in terms of $\bar{\chi}(u)$ is
not always an increasing function of the correlation parameter.
The special cases in these figures, in light of Proposition \ref{pro1}, are the ones with
$0\leq\delta_2<\delta_1$ and $\delta_2<0\leq\delta_1$
which indicates some evidence for asymptotic dependence.
Lastly, they manifestly exhibit positive extremal dependence while
the negative and near dependence is appeared to a lesser extent
under the current setting. Overall, these findings indicate that this version of
GSN distribution can induce tail
dependence at sub-asymptotic and/or asymptotic
levels.
The strength of dependence appears to vary with correlation as well as
separation of skewness parameters. Of course, the parameters value
$\delta_1=0$ and $\delta_2=0$ display null extremal dependence.
Notably, the tail probabilities
of all pairs of variables under the GSN distribution is not driven
by a single parameter like a multivariate t distribution whose UTDC
converges to the positive value for finite value of the degrees of
freedom parameter regardless of the correlation parameter \citep{Dem}.

The numerical integration is implemented using
the statistical software R \citep{Rcor}. The
computation of the $\bar{\chi}(\cdot)$ is restricted to the interval
$[\zeta,1-\zeta]$, where $\zeta$ is defined by machine precision,
about $10^{-9}$ on our machine.
\section{The GSN distributed random field}\label{model2}
A spatial RF model $\{Y(\bm s):\bm s\in
\mathcal{D}\subset\mathbb{R}^2\}$ is a collection of random
variables indexed via $\bm s\in \mathcal{D}$. If a finite set of
locations $\{\bm s_1,\ldots,\bm s_n\}\in\mathcal{D}$ is observed,
then the finite-dimensional distributions for each $n\geq1$ must
satisfy Kolmogorov's compatibility conditions.
Using the mgf in \eqref{mgfssn}, it can be shown that the
compatibility conditions are satisfied under the GSN finite-dimensional
distributions. Therefore, the GSN distributed RF is well-defined
\citep{Mah2}.

A key idea to link the GSN distribution and the spatial model is to
view $Y$ as follows
\begin{align}
Y(\bm s)=\mu(\bm s)+Z(\bm s)+\epsilon(\bm s),\label{SLM}
\end{align}
where $\mu(\bm s)$, $Z(\bm s)$ and $\epsilon(\bm s)$ being location
dependent mean function, a smooth-scale SGRF with
\begin{align*}
\GSN_n(\bm0,\bm\Sigma,\bm\delta)
\end{align*}
finite-dimensional distribution and an i.i.d. (independent and identically distributed)
GRF independent of $Z(\bm s)$ with mean zero and variance
$\tau^2$, respectively. The parameter $\tau^2$ is called nugget
effect in geostatistical context. Equivalently, let that $\{W(\bm
s): \bm s\in\mathbb{R}^2\}$ and $\{\delta(\bm s): \bm s\in\mathbb{
R}^2\}$ be two stationary GRF defined on $\mathbb{R}^2$ and assume
that $W(\bm s)$ and $\delta(\bm s)$ be independent, with components
having following means
\begin{align*}
\E[W(\bm s)]=0,\quad \E[\delta(\bm s)]=\gamma, \quad \bm
s\in\mathbb{R}^2,
\end{align*}
and covariances
\begin{align}
\Cov[W(\bm s),W(\bm s')]=\sigma^2\rho_w(||\bm s-\bm s'||),\quad
\Cov[\delta(\bm s),\delta(\bm s')]=\gamma^2\rho_\delta(||\bm s-\bm
s'||),\quad \bm s,\bm s'\in \mathbb{R}^2,\label{covDW}
\end{align}
where $||\bm s-\bm s'||$ is a Euclidean distance between two field
measurement locations $\bm s$ and $\bm s'$. Furthermore,
$\rho_w(\cdot)$ and $\rho_\delta(\cdot)$ in \eqref{covDW} are the
corresponding spatial correlation functions of $W(\bm s)$ and
$\delta(\bm s)$, respectively. Note that, under this formulation when
distance between spatial locations goes to the infinity, the
covariances of this SGRF vanish. One option for the
correlation function is the Mat$\acute{\textrm{e}}$rn family of the
correlation functions, given by
\begin{align}
\rho(||\bm s-\bm
s'||;\psi,\xi)=\frac{1}{2^{\xi-1}\Gamma(\xi)}\left(\frac{1}{\psi}||\bm
s-\bm s'||\right)^{\xi}\mathcal{K}_\xi\left(\frac{1}{\psi}||\bm
s-\bm s'||\right),\label{MAT}
\end{align}
where $\psi$ is a range parameter, $\Gamma(\cdot)$ is the gamma
function and $\mathcal{K}_\xi(\cdot)$ is the modified Bessel
function of the third kind of order $\xi>0$ \citep{St}. It depends
on a smoothness parameter $\xi$ which directly controls the mean
square differentiability of RF realizations. If $\xi>1$ then
Mat$\acute{\textrm{e}}$rn correlation functions are once mean square
differentiable, and if $\xi=3/2$, the correlation functions are of
the closed form $(1+||\bm s-\bm s'||/\psi)\exp(-||\bm s-\bm
s'||/\psi)$. Because parameters of the Mat$\acute{\textrm{e}}$rn
correlation function often being poorly
identified \citep{Zha}, we set $\xi=3/2$. Under this setting, the
range parameter $\psi$ is appeared out of the exponential term.
Accordingly, data may contribute more information on estimation of
$\psi$. Throughout the text, similar Mat$\acute{\textrm{e}}$rn correlation
function with $\xi=3/2$ is adopted for each of $\rho_w(\cdot)$ and
$\rho_\delta(\cdot)$ to ensure the existence of such GSN distributed RF model.

By assuming $T(\bm s)=V(\bm s)-\sqrt{2/\pi}$ and
$V(\bm s)\stackrel{\;\rm i.i.d.}{\sim}\HN_1(0,1)$, one can employ the
following representation of the spatial SGRF model
\begin{align}
Y(\bm s)=\mu(\bm s)+W(\bm s)+\delta(\bm s)T(\bm s)+\epsilon(\bm
s).\label{SSNM}
\end{align}

Several properties of SGRF in the GSN family could be deduced from
\eqref{SSNM}. Let us assume the case when the mean function is
constant: $\mu(\bm s)=\mu$ for all $\bm s\in\mathbb{R}^2$, then the
skewed field defined in \eqref{SSNM} is stationary with expectations
\begin{align*}
\E[Y(\bm s)]=\mu,\quad \bm s\in\mathbb{R}^2,
\end{align*}
variances
\begin{align*}
\Var[Y(\bm
s)]=\tau^2+\sigma^2+2\gamma^2\left(1-\frac{2}{\pi}\right),\quad \bm
s\in\mathbb{R}^2,
\end{align*}
and covariances
\begin{align*}
\Cov[Y(\bm s),Y(\bm s')]=\sigma^2{\rm \rho}_w(||\bm s-\bm
s'||),\quad \bm s,\bm s'\in\mathbb{R}^2.
\end{align*}
\section{Heavy-tailed construction of SGRF model}\label{model3}
In this section, we consider the scale mixture of the
GSN distributed RF model and discuss about its skewness and kurtosis.
Let $Y_i=Y(\bm s_i)$ denote the observation at spatial location
$\bm s_i$ and consider the data model
\begin{align}
Y_i=\bm X(\bm
s_i)'\bm\beta+\sigma\lambda_i^{-1/2}W_i+\gamma\lambda_i^{-1/2}\delta_i
T_i+\epsilon_i,\label{SGLG}
\end{align}
where $X(\bm s_i)$ is a $(p\times1)$ vector of known location
dependent covariates, $\bm\beta$ is a $(p\times1)$ vector of unknown
regression parameters, $\sigma>0$ is the scale parameter,
$\gamma\in\mathbb{R}$ is the asymmetry parameter and
the spatial random effect $W_i=W(\bm s_i)$ is independent of
$T_i=T(\bm s_i)$, follows a GRF model defined by
\begin{align*}
\bm W=(W_1,\ldots,W_n)'\sim \N_n(0,\bm H),
\end{align*}
where $\bm H$ is a correlation matrix describing the spatial
dependence, whose elements are given by \eqref{MAT} with $\xi=3/2$. The asymmetry
parameter, $\gamma$, controls the direction and magnitude of
skewness. It is worth recalling that negative values of $\gamma$
induce negative skewness, positive values generate positive
skewness, and $\gamma=0$ corresponds to symmetry. To allow altered
skewness mixing variable for each measurement location, by setting
$\delta_i=\delta(\bm s_i)$, we adopt the following GRF model for
$\delta(\cdot)$:
\begin{align*}
\bm\delta=(\delta_1,\ldots,\delta_n)'\sim\N_n(\bm 1_n,\bm H).
\end{align*}

Also an additional set of latent random variables,
$\lambda_i=\lambda(\bm s_i)>0$, is introduced here to deal with the presence of
outliers in the spatial responses. The scale mixing variables,
$\bm\lambda=(\lambda_1,\ldots,\lambda_n)'$, are assumed to be
spatially correlated to induce the mean square continuity. Hence,
$\bm\lambda$ on the logarithm scale can be modeled by a GRF model as
follows
\begin{align}
\ln(\bm\lambda)\sim\N_n\left(-\frac{\nu}{2}\bm 1_n,\nu\bm
H\right),\label{GLG}
\end{align}
where $\nu>0$ is the tail-weight parameter regulating the heaviness
in tails of the RF. When the spatial model is rescaled according to
the transformed GRF model in \eqref{GLG}, on average the resultant
finite-dimensional distribution is a GSN probability model with
inflated variance. While large values of the tail-weight parameter
has been found to provide heavier tails, the SGRF model corresponds
to the limiting case when $\nu$ tends toward zero. Additionally,
taking into account different correlation structures
entails little flexibility in comparison to the skewness
parameters. Therefore, we assume similar spatial correlation matrix
for $\lambda(\cdot)$ and $W(\cdot)$.
This assumption further can be assessed by analyzing real data sets
in terms of better prediction results.
\begin{figure*}[!t]
\vspace{0mm}
\begin{center}
\includegraphics[scale=0.8,angle=0]{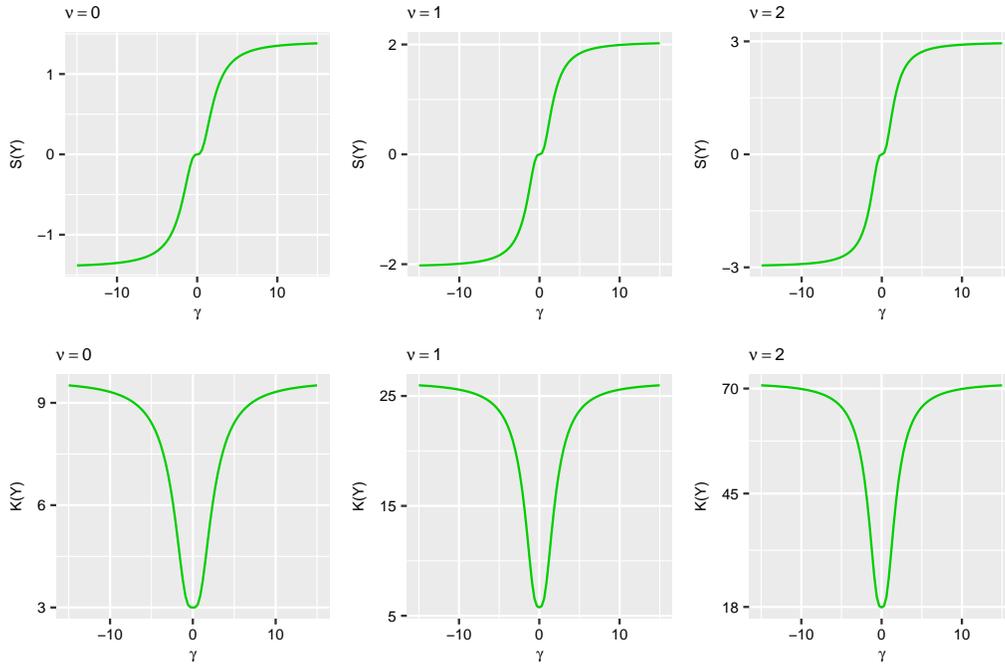}
 \end{center}
 \caption{Coefficients of skewness (first row) and kurtosis (second row) for the model in (\ref{SGLG}) with respect to various values of the parameters
 $\gamma$ and $\nu$.}\label{sk-ku}
\end{figure*}

Now, we compute the coefficient of skewness and kurtosis of
\eqref{SGLG}. Tedious but straightforward algebra proved that these
coefficients are given, respectively, by
\begin{align*}
{\rm
S}(Y)&=\frac{4\gamma^3\exp\left(\frac{15\nu}{8}\right)\left(\frac{4}{\pi}-1\right)\sqrt{\frac{2}{\pi}}}
{\left[\tau^2+\sigma^2\exp(\nu)+2\gamma^2\exp(\nu)\left(1-\frac{2}{\pi}\right)\right]^{3/2}},\\\\
{\rm
K}(Y)&=\frac{3\tau^4+3\sigma^4\exp(3\nu)+6\sigma^2\tau^2\exp(\nu)+12\sigma^2\gamma^2A_1+
12\tau^2\gamma^2A_2+10\gamma^4A_3}
{\left[\tau^2+\sigma^2\exp(\nu)+2\gamma^2\exp(\nu)\left(1-\frac{2}{\pi}\right)\right]^2},
\end{align*}
in which
\begin{align*}
A_1= \exp(3\nu)\left(1-\frac{2}{\pi}\right),\quad
A_2=\exp(\nu)\left(1-\frac{2}{\pi}\right),\quad
A_3=\exp(3\nu)\left(3-\frac{4}{\pi}-\frac{12}{\pi^2}\right).
\end{align*}
These measures are used to produce \eqref{sk-ku} which is
plotted by varying $\gamma$ for different values of the parameter
$\nu$ after assuming $\tau=1$ and $\sigma=1$, denoted by standard case. We do
not introduce these parameters into the assessments, because they have
minor impacts on skewness and kurtosis. Essentially, for fixed values
of $\gamma$ and $\nu$, ${\rm S}(Y)$, and ${\rm K}(Y)$ decrease for
increasing values of $\tau$ or $\sigma$. However, for large values
of $\tau$ or $\sigma$, the skewness and kurtosis coefficients have
similar patterns for greater values of $\gamma$ and $\nu$ than
the ones that were assumed in the standard case. It can be
examined when $\nu=0$, ${\rm S}(Y)$ varies in the interval
$(-1.41,1.41)$, while ${\rm K}(Y)$ takes values in the range
$(3,9.67)$. An aspect to be stressed right away is that the
parameter $\nu$ ($\gamma$) is not concerned with tail (skewness) of
the RF model merely. The addition of this extra parameter, $\nu$
($\gamma$), to allow for flexibility in the GRF supports the
skewness (tail) whose size expands with $\nu$ ($\gamma$).
Consequently, both of $\nu$ and $\gamma$ control the non-Gaussian
strength.
\section{Conclusions}\label{conclu}
We extend the methodology previously presented in the
literature to accurately incorporate size and direction
of the skewness into the RF model. An elegant consequence
is that the SGRF model with spatially varying skewness parameters do
a great job in capturing tail probabilities which are thought to
possess some degree of tail dependence.

 \end{document}